# Approximation, regularization and smoothing of trigonometric splines

Denysiuk V.P. Dr of Phys-Math. sciences, Professor, Kiev, Ukraine
National Aviation University
kvomden@nau.edu.ua

## Annotation

The methods of approximation, regularization and smoothing of trigonometric interpolation splines are considered in the paper. It is shown that trigonometric splines can be considered from two points of view - as a trigonometric Fourier series and as discrete trigonometric Fourier series according to certain systems of functions that are smoothness carriers. It is argued that with approximation and smoothing of trigonometric splines it is expedient to consider as discrete rows, since their differential properties are stored.

**Keywords:** Trigonometric splines, continuum and discrete Fourier series approximation, regularization, smoothing.

## Introduction

The problem of approximation of functions arises when solving many practical and theoretical problems. Thus, the approximation of functions is an important auxiliary apparatus in solving other problems of numerical analysis - numerical integration and differentiation, finding approximate solutions of differential and integral equations, and so on. The process of approaching one function to another is often called an approximation.

There are continuous and discrete approach to the approximation problem. The continuous approximation is applicable to functions, the representation of which is known throughout the function interval. The discrete approximation is applied to the functions given by the values in nodes of some, as a rule, uniform grids set on an interval of task of function.

Since trigonometric splines belong to the class of periodic functions, in what follows we will confine ourselves to considering the problem of approximating only periodic functions with period $2\pi$.

The task of continuing approximation of some continuous, periodic with period $2\pi$ function $f(x)$ can be formulated as follows. Find the trigonometric polynomial

$$T_m(x) = \frac{a_0}{2} + \sum_{k=1}^{m}\left(a_k \cos kx + b_k \sin kx\right)$$

Of an order $m$ ($m = 1, 2, ...$), on which the smallest of the values is achieved

$$\delta(m) = \int_{0}^{2\pi}\left[T_m(x) - f(x)\right]^2 dx.$$

Of course, we assume that the function $f(x)$ integrated with square.

From the course of mathematical analysis it is known that the smallest value $\delta(m, p_0, p_1, ..., p_{2m})$ is achieved in the case when the polynomial $T_m(x)$ is the partial sum of the trigonometric Fourier series, that is its coefficients $a_0, a_1, b_1, ..., a_k, b_k, ...$ are the coefficients of the Fourier function $f(x)$, which are calculated by formulas

$$a_0 = \frac{1}{\pi}\int_{0}^{2\pi} f(x)dx;$$

$$a_k = \frac{1}{\pi}\int_{0}^{2\pi} f(x)\cos kx\, dx, \qquad b_k = \frac{1}{\pi}\int_{0}^{2\pi} f(x)\sin kx\, dx, \quad (k = 1, 2, ...)$$

The problem of discrete approximation of periodic functions is often considered in such a statement [1]. Let on the segment $[0, 2\pi)$ the function is set $f(x)$ of a certain class. Suppose that a uniform grid is also given on this segment $\Delta_N = \{x_j\}_{j=1}^N$, $x_j = \frac{2\pi}{N}(j-1)$, $N = 2n+1$, ($n = 1, 2, ...$) and are known values $f(x_j) = f_j$ ($j = 1, 2, ..., N$) this function in the grid nodes $\Delta_N$. It is necessary to find a trigonometric polynomial

$$T_m^*(x) = \frac{a_0^*}{2} + \sum_{k=1}^m \left(a_k^* \cos kx + b_k^* \sin kx\right)$$

Of an order $m$ ($m \leq n$), on which the smallest value is reached

$$\delta^*(m) = \sum_{j=1}^m \left[T_m^*(x_j) - f_j\right]^2.$$

Note that the method by which the value is minimized $\delta^*(m)$ often called the least squares method.

Together with the continuum theory of trigonometric Fourier series, there is a discrete theory of such series [1], [2], from which it follows that the smallest value of the quantity $\delta^*(m)$ is achieved in the case where the coefficients $a_0^*, a_1^*, b_1^*, ..., a_m^*, b_m^*$ of polynomial $T_m^*(x)$ are discrete Fourier coefficients calculated by Bessel's formulas [1]

$$a_0^* = \frac{2}{N} \sum_{j=1}^N f_j \ ;$$

$$a_k^* = \frac{2}{N} \sum_{j=1}^N f_j \cos kt_j, \quad b_k^* = \frac{2}{N} \sum_{j=1}^N f_j \sin kt_j,$$

$$k = 1, 2, ..., m, \quad m \leq n.$$

Note that these formulas can be obtained from the known formulas for the Fourier coefficients of the function $f(x)$, if we calculate the integrals included in them using composite formulas of rectangles [1] or trapezoids; also these formulas can be obtained as coefficients of the interpolation trigonometric polynomial

It is clear that the value $\delta^*(m)$ significantly depends on the value of the parameter $m$. So, when $m = n$ we can determine these coefficients so that the conditions are met

$$T_n^*(x_j) = f_j, \quad j = 1, 2, ..., N \ ;$$

and $\delta^*(n) = 0$. This approximation problem is called the interpolation problem, the criterion of equality of polynomial values by the value of the function $f(x)$ in grid nodes $\Delta_N$ called the interpolation criteria and the polynomial $T_n^*(x)$, that satisfies these conditions is called an interpolation polynomial. The theory of construction of such polynomials is well developed and is often applied in practice (see, for example, [1,2]).

The application of the interpolation criterion to construct the approximation assumes that the value of the function $f_j$, ($j = 1, 2, ..., N$), known for sure. However, in many mathematical modelling problems, these values of the function are obtained as a result of physical measurements and are not accurate; in such cases, apply the method of least squares and consider the value $\delta^*(m)$, ($m < n$). Hereinafter, the least squares approximation will be called simply an approximation.

# The purpose of the work

Development of methods for approximation, regularization and smoothing of trigonometric splines

# The main part

### Approximation of trigonometric splines

Previously entered grid $\Delta_N = \{x_j\}_{j=1}^N$ in the future we will denote as $\Delta_N^{(0)} = \{x_j^{(0)}\}_{j=1}^N$. We will also introduce a grid $\Delta_N^{(1)} = \{x_j^{(1)}\}_{j=1}^N$, $x_j^{(1)} = \frac{\pi}{N}(2j-1)$.

In [3] it was shown that in the general case trigonometric interpolation splines on grids $\Delta_N^{(I)}$ can be submitted in this way

$$St^{(I_1,I_2)}(\Gamma,H,\nu,r,x) = \frac{a_0^{*(I_2)}}{2} + \sum_{k=1}^{\frac{N-1}{2}}\left[a_k^{*(I_2)}\frac{C_k^{(I_1)}(\Gamma,\nu,r,x)}{hc^{(I_1,I_2)}(\Gamma,r,k)} + b_k^{*(I_2)}\frac{S_k^{(I_1)}(H,\nu,r,x)}{hs^{(I_1,I_2)}(H,r,k)}\right],$$

where

$$C_k^{(I_1)}(\Gamma,\nu,r,x) =$$
$$= \gamma_1\nu_k(r)\cos kx + \sum_{m=1}^{\infty}(-1)^{mI_1}\left[\gamma_3\nu_{mN+k}(r)\cos(mN+k)x + \gamma_2\nu_{mN-k}(r)\cos(mN-k)x\right];$$

$$S_k^{(I_1)}(H,\nu,r,x) =$$
$$= \eta_1\nu_k(r)\sin kx + \sum_{m=1}^{\infty}(-1)^{mI_1}\left[\eta_3\nu_{mN+k}(r)\sin(mN+k)x - \eta_2\nu_{mN-k}(r)\sin(mN-k)x\right],$$

with convergence factors $\nu_k(r)$, having a descending order $O(k^{-(1+r)})$, and interpolation factors

$$hc_k^{(I_1,I_2)}(\Gamma,\nu,r) = \gamma_1\nu_k(r) + \sum_{m=1}^{\infty}(-1)^{m(I_1-I_2)}\left[\gamma_3\nu_{mN+k}(r) + \gamma_2\nu_{mN-k}(r)\right],$$

$$hs_k^{(I_1,I_2)}(H,\nu,r) = \eta_1\nu_k(r) + \sum_{m=1}^{\infty}(-1)^{m(I_1-I_2)}\left[\eta_3\nu_{mN+k}(r) + \eta_2\nu_{mN-k}(r)\right].$$

where the indicator $I_1$ ($I_1 = 0,1$) determines the stitching the grid, indicator $I_2$ ($I_2 = 0,1$) determines the interpolation grid, $a_0^{*(I_2)}, a_k^{*(I_2)}, b_k^{*(I_2)}$ - The coefficients of interpolation trigonometric polynomials on the grid $\Delta_N^{(I_2)}$, $r$, ($r = 1,2,...$) - parameter that determines the order of the spline, and $\Gamma = \{\gamma_1,\gamma_2,\gamma_3\}$ and $H = \{\eta_1,\eta_2,\eta_3\}$ - parameter vectors, and parameters $\gamma_k$ and $\eta_k$, ($k = 1,2,3$) take non-zero real values. Note that to reduce the notation in the splines and the functions through which they are built, we omit the dependence on the number $N$ nodes of interpolation grids $\Delta_N^{(I)}$.

In this paper, without losing generality, we will limit ourselves to the consideration of the trigonometric spline, which is sewn in the nodes of the grid $\Delta_N^{(0)}$ and interpolates the function $f(x)$ in nodes of the same grid; such a spline can be presented in a simpler form

$$St(\Gamma,H,\nu,r,x) = \frac{a_0}{2} + \sum_{k=1}^{\frac{N-1}{2}}\left[a_k\frac{C_k(\Gamma,\nu,r,x)}{hc(\Gamma,r,k)} + b_k\frac{S_k(H,\nu,r,x)}{hs(H,r,k)}\right],$$

where

$$C_k(\Gamma,\nu,r,x) =$$
$$= \gamma_1\nu_k(r)\cos kx + \sum_{m=1}^{\infty}\left[\gamma_3\nu_{mN+k}(r)\cos(mN+k)x + \gamma_2\nu_{mN-k}(r)\cos(mN-k)x\right];$$

$$S_k(H,\nu,r,x) =$$
$$= \eta_1\nu_k(r)\sin kx + \sum_{m=1}^{\infty}\left[\eta_3\nu_{mN+k}(r)\sin(mN+k)x - \eta_2\nu_{mN-k}(r)\sin(mN-k)x\right],$$

with factors of convergence $v_k(r)$, having a descending order $O(k^{-(1+r)})$, and interpolation factors

$$hc_k(\Gamma,v,r) = \gamma_1 v_k(r) + \sum_{m=1}^{\infty}\left[\gamma_3 v_{mN+k}(r) + \gamma_2 v_{mN-k}(r)\right],$$

$$hs_k(H,v,r) = \eta_1 v_k(r) + \sum_{m=1}^{\infty}\left[\eta_3 v_{mN+k}(r) + \eta_2 v_{mN-k}(r)\right].$$

where $a_0, a_\kappa, b_\kappa$ - coefficients of the interpolation trigonometric polynomial on the grid $\Delta_N^{(0)}$, $r$, ($r = 1, 2, ...$) - parameter that determines the order of the spline, and $\Gamma = \{\gamma_1,\gamma_2,\gamma_3\}$ and $H = \{\eta_1,\eta_2,\eta_3\}$ - the same as before the parameter vectors.

It is clear that trigonometric splines have a double structure; On the one hand, they are a function whose representation in the form of an infinite Fourier series is known over the entire interval of the approximate function. On the other hand, such splines can be considered as finite amount by function $\dfrac{C_k(\Gamma,v,r,t_j)}{hc(\Gamma,r,k)}$ and $\dfrac{S_k(H,v,r,t_j)}{hs(H,r,k)}$, For which at all values $\Gamma, H, v, r$ equations are fulfilled

$$\frac{C_k(\Gamma,v,r,t_j)}{hc(\Gamma,r,k)} = \cos kt_j; \quad \frac{S_k(H,v,r,t_j)}{hs(H,r,k)} = \sin kt_j, \quad (j=1,2,\ldots,N).$$

It follows that trigonometric splines can also be considered as sums of a discrete Fourier series.

Thus, both continuous and discrete approximations can be applied to trigonometric splines. Consider the problem of approximation of trigonometric splines in more detail.

First of all, note that features $C_k(\Gamma,v,r,t)$ and $S_k(H,v,r,t)$, included in the spline representation $St(\Gamma,H,v,r,t)$, are in fact uniformly convergent Fourier series and are carriers of the differential properties of trigonometric splines. It follows that by applying the classical Fourier series approximation, which consists in discarding the residues of these series, we lose the differential properties of trigonometric splines. Therefore, in those problems where these properties must be preserved, the classical approximation of trigonometric series is impractical.

Let us now consider the discrete approximation of trigonometric splines by the least squares method. Because the functions $\dfrac{C_k(\Gamma,v,r,t_j)}{hc(\Gamma,r,k)}$ and $\dfrac{S_k(H,v,r,t_j)}{hs(H,r,k)}$ form an orthogonal in the discrete sense on the grid system of functions, then, according to the general theory, the approximation by the method of least squares is to discard the higher coefficients $a_0^*, a_1^*, b_1^*, \ldots, a_m^*, b_m^*$ trigonometric spline. It is clear that such approximation preserves the differential properties of the trigonometric spline.

## Approximation of functions with additional constraints

Applications often require finding a solution to the problem of approximating of a function $f$ with additional restrictions imposed on the approximating function $g$. For example, in the classical formulation of the problem of approximation of a function $f$ trigonometric polynomials, additional restrictions are the requirements that the degree of the approximating polynomial does not exceed some parameter $m$, ($m = 1, 2, ...$). Additional limitations of another type are the requirements for the approximation function $g$ had continuous derivatives of a certain order (see [4], [5]), or the requirement that the approximating function $g$ coincided with the function $f$ not only on the set of given discrete points, but also in some other sense. Consider the problem of approximating functions with additional constraints in more detail. In this case, in the role of approximating functions, we will consider only the classes of trigonometric functions.

Additional limitations that are required to have an approximating function $g$ had a fairly smooth derivative of the order $p$, ($p = 1, 2, ...$), can be reduced to the task of minimizing the functionality

$$\Psi(g) = \int_0^{2\pi}\left[g^{(p)}(x)\right]^2 dx$$

In this case, the solution of the problem of approximation of the function $f$ it is expedient to look in the form of function $g$, on which the minimum of functionality is reached

$$\Phi_\lambda^p(f,g) = \int_0^{2\pi} \left\{ [f(x) - g(x)]^2 + \lambda^{2p} \left[ g^{(p)}(x) \right]^2 \right\} dx$$

where $\lambda, p, (n = 1, 2, ...)$ - some parameters.

The problem of function approximation $f$, at which it is required to minimize functionality $\Phi_\lambda^p(f,g)$ called the regularization problem, and the parameters $\lambda$ and $p$ called regularization parameters [6].

Applying the methods of variational calculus, we can make a differential Euler equation with respect to the desired function $g$, which has the form

$$\lambda^{2p} g^{(2p)} + (-1)^p (g - f) = 0$$

If the function $f$ can be represented as a Fourier series,

$$f(x) = \frac{a_0}{2} + \sum_{k=1}^{\infty} [a_k \cos kx + b_k \sin kx]$$

then a direct check can make sure [6] that the function $g$, which is the solution of the regularization problem can be represented as

$$g(x) = \frac{a_0}{2} + \sum_{k=1}^{\infty} \frac{1}{\left[1 + (\lambda k)^{2p}\right]} [a_k \cos kx + b_k \sin kx] /$$

In other words, the Fourier coefficients of the function $g(x)$ which is the solution of the regularization problem, are obtained from the Fourier coefficients of the function $f(x)$ by multiplying them by the values of

$$\tau_k(\lambda, p) = \frac{1}{1 + (\lambda k)^{2p}}, \quad k = 1, 2, ... \ .$$

Graph of values $\tau_k(\lambda, p)$ is shown on Fig.1.

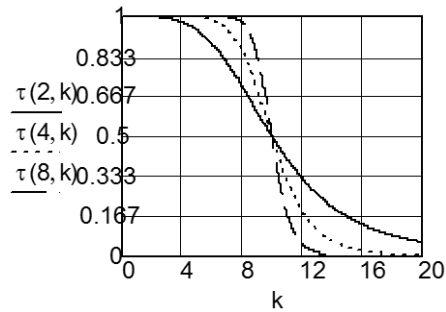

Fig.1. Graph of values $\tau_k(\lambda, p)$ at $\lambda = .1$

This result is easy to transfer to trigonometric splines. It is clear that in the role of the function $f(x)$ we consider functions $C_k(\Gamma, \nu, r, t)$ and $S_k(H, \nu, r, t)$, which take the form:

$$C_k R_{\lambda,p}(\Gamma, \nu, r, x) = \frac{1}{1 + (\lambda k)^{2p}} \gamma_1 \nu_k(r) \cos kx +$$

$$+ \sum_{m=1}^{\infty} \left[ \frac{1}{1 + [\lambda(mN+k)]^{2p}} \gamma_3 \nu_{mN+k}(r) \cos(mN+k)x + \frac{1}{1 + [\lambda(mN-k)]^{2p}} \gamma_2 \nu_{mN-k}(r) \cos(mN-k)x \right];$$

$$S_k R_{\lambda,p}(H, \nu, r, x) = \frac{1}{1 + (\lambda k)^{2p}} \eta_1 \nu_k(r) \sin kx +$$

$$+ \sum_{m=1}^{\infty} \left[ \frac{1}{1 + [\lambda(mN+k)]^{2p}} \eta_3 \nu_{mN+k}(r) \sin(mN+k)x - \frac{1}{1 + [\lambda(mN-k)]^{2p}} \eta_2 \nu_{mN-k}(r) \sin(mN-k)x \right].$$

When applying the regularization of trigonometric splines, it is necessary to consider that the factors $\tau_k(\lambda,p) = \dfrac{1}{1+(\lambda k)^{2p}}$ will significantly increase the order of decrease the Fourier coefficients of these functions, and therefore change their differential properties. It follows that the regulatory task is appropriate to use for splines with relatively small values of the parameter $r$. In addition, the important role is played by the parameter $\lambda$, which determines the coefficient number that is weakened twice.

## Smoothing trigonometric splines

In the study of finite sums of Fourier series, linear methods are often used - the summation of these amounts [7], which are variants of the Chesaro's method of generalized summation of rows; These methods significantly improve the properties of these amounts [9]. Since trigonometric splines, as we have already said above, are the finite sum of discrete Fourier series, then there is an opportunity to use the methods of summing up to these splines. Trigonometric splines with introduced in them $\lambda_k^{(n)}$ - Multipliers will be called smoothed splines, indicating the name of multipliers. So, for example, trigonometric splines with entered in them $\lambda_k^{(n)}$ - The factors of the Fejér will be called splines smoothed by Fejér.

Trigonometric spline, smoothed with $\lambda_k^{(n)}$ - multipliers can be submitted in the form

$$St(\Gamma, H, \nu, r, \Lambda, x) = \frac{a_0^*}{2} + \sum_{k=1}^{\frac{N-1}{2}} \lambda_k^{(n)} \left[ a_k^* \frac{C_k(\Gamma,\nu,r,x)}{hc(\Gamma,r,k)} + b_k^* \frac{S_k(H,\nu,r,x)}{hs(H,r,k)} \right].$$

Known $\lambda_k^{(n)}$ Multipliers of the Fejér, Vallée-Poussin, Hemming-Tyuki, Bohman, Riman-Lanczos, Rogosinski, etc. [8]. We have used a modified Fejér's factor that looks like

$$\lambda_k(\alpha,n) = \left[1 - \frac{k}{n+1}\right]^\alpha, \quad (\alpha > 0),$$

and which is rather flexible depending on the values of the parameter $\alpha$. The graph of such a factor is shown in Fig.2

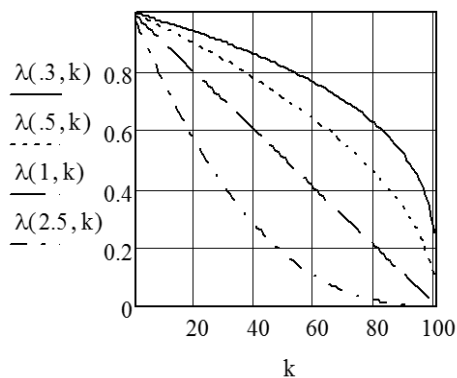

Fig. 2. A graph of a modified factor of the Fejér at $n = 100$.

## Smoothing data

Finally, another method of smoothing trigonometric splines is a method in which not the spline is smoothed, but the value of an interpolated function; Such a smoothing method is called data smoothing. This method is considered in detail in [1], [2] and there are also computing formulas for data smoothing. In [10] it is shown that data smoothing can be interpreted as one of the variants of the data filtering process using digital filters, in which the distribution properties of spline are not changed.

Finally, we note that the methods of approximation, regularization and smoothing of trigonometric splines can be recommended for widespread use in practice.

## Conclusions

1. Approximation, regularization and smoothing of trigonometric splines are considered.
2. It has been shown that trigonometric splines can be considered from two points - as a trigonometric series of Fourier and as the sum of discrete trigonometric Fourier series.
3. With approximation of trigonometric splines, it is advisable to consider them as the sum of discrete Fourier series, as their differential properties are preserved.
4. Regularization of trigonometric splines changes them differential properties; This should be taken into account for further studies of these splines.
5. When smoothing trigonometric splines by summation methods it is expedient to consider them as the sums of discrete Fourier series.
6. In many tasks associated with a generalized summation of trigonometric ranks, it is advisable to consider the modified Fejér's multiplier.
7. Of course, the proposed methods of approximation, regularization and smoothing of trigonometric splines require further research.

## List of references